\documentclass[12pt]{article}
\usepackage[T2A]{fontenc}
\usepackage[utf8]{inputenc}
\usepackage{amssymb,amsmath,amsfonts,amsthm,amscd,latexsym,verbatim,graphics,epsfig,indentfirst,xcolor}
\usepackage{geometry}
\def\rb{\mathrm{rb}}

\def\charr{\mathrm{char}\,}
\def\rank{\mathrm{rank}}

\def\Aut{\mathrm{Aut}}

\def\Imm{\mathrm{Im}}

\begin{document}

\sloppy

\begin{center}
{\Large
Rota---Baxter operators of weight zero \\
on upper-triangular matrices of order three}

V. Gubarev
\end{center}

\begin{abstract}
We describe all Rota---Baxter operators~$R$ of weight zero
on the algebra~$U_3(F)$ of upper-triangular matrices of order three over a field of characteristic~0.
For this, we apply the following three ingredients: properties of~$R(1)$, conjugation with suitable (anti)automorphisms of~$U_3(F)$, and computation with the help of \texttt{Singular} of the Gr\"{o}bner basis for the system of the equations arisen on the coefficients of~$R$.

\end{abstract}

\section{Introduction}

Rota---Baxter operators were defined by G. Baxter in 1960~\cite{Baxter}.
In the 1960--1970s, the main interest to such operators was concerned to Banach commutative and associative algebras.
In the 1980s, the connection between classical and modified Yang---Baxter equations and Rota---Baxter operators on semisimple Lie algebras was found~\cite{BelaDrin82,Semenov83}.
Since the 2000s, the active study of Rota---Baxter operators on associative
algebras (without topology restrictions) has begun; see the monograph~\cite{GuoMonograph} written by Li Guo.

One of the interesting directions in the study of Rota---Baxter operators is the problem of their classification on a given algebra. RB-operators of any weight were classified on 2-dimensional associative algebras over~$\mathbb{C}$~\cite{preLie}, of nonzero weight all and of
weight zero some of them were classified on 3-dimensional associative algebras over~$\mathbb{C}$~\cite{Abdujabborov,AnBai}.

Description of Rota---Baxter operators of any weight on $M_2(\mathbb{C})$ (see~\cite{BGP} and references therein),
of nonzero weight on $M_3(\mathbb{C})$ in the series of three works~\cite{GonGub,Gub2021,Gub2024} was obtained.
In 2013, V.V.~Sokolov described all skew-symmetric RB-operators of weight zero on $M_3(\mathbb{C})$~\cite{Sokolov}.
The general interest to Rota---Baxter operators on matrix algebras is clarified by their deep connection with the associative Yang---Baxter equation~\cite{Aguiar00} and double Poisson algebras~\cite{DoubleLie,Schedler}.

In~\cite{Mat2}, the approach involving Gr\"{o}bner bases to deal with Rota---Baxter operators was successfully applied to state the description of RB-operators of weight zero on~$M_2(\mathbb{C})$.
However, if one tries to use it for algebras of higher dimensions, computational problems arise. 
It means that the obtained system of quadratic relations may be huge, and it is difficult to solve it.

In the current work, we suggest a new approach to classify Rota---Baxter operators on a given finite-dimensional unital associative algebra~$A$. The spectral properties of an RB-operator~$R$ on~$A$ as well as the properties of~$R(1)$ give important restrictions on the possible form of~$R$~\cite{Unital}.
If the weight of~$R$ is zero, we have that $R(1)$ is a nilpotent element of~$A$, and $R$ is a~nilpotent operator on~$A$.
Thus, we may split the general classification into subcases depending on the form of~$R(1)$ and only then apply Gr\"{o}bner bases technique.

We study the algebra~$U_3(F)$ of upper-triangular matrices of order three.
Surprisingly, the stated above strategy works well, and the classification of RB-operators of weight zero on~$U_3(F)$ is fulfilled.
We get 40 cases, except for the two operators
\begin{gather*}
R(e_{12}) = R(e_{13}) = 0, \quad
R(e_{11}),
R(e_{22}),
R(e_{33}),
R(e_{23})\in L(e_{12},e_{13}); \\
R(e_{13}) = 0, \quad
R(e_{11}), R(e_{12}), R(e_{22}), R(e_{23}), R(e_{33})\in L(e_{13}),
\end{gather*}
which are Rota---Baxter ones from the known arguments,
all other 38 cases involve maximum two parameters.
Let us highlight that it is important to conjugate RB-operators by an (anti)automorphism of~$U_3(F)$.
We apply (anti)automorphisms of~$U_3(F)$ twice: for the choice of a nilpotent matrix $R(1)$ and afterwards,
analyzing the solutions of obtained by \texttt{Singular}~\cite{DGPS} systems of equations.

\section{Preliminaries}

Let $F$ stand for the ground field.

We give a list of some known properties of Rota---Baxter operators.

{\bf Lemma 1}~\cite{GuoMonograph}.
Given an RB-operator $P$ of weight~0 and $k\in F\setminus\{0\}$,
the operator $k^{-1}P$ is again an RB-operator of weight~0.

Below, by an antiautomorphism of $A$ we mean a bijection~$\psi$ from $A$ to $A$ satisfying $\psi(xy) = \psi(y)\psi(x)$ for all $x,y\in A$;
e.g., transpose on the matrix algebra.

{\bf Lemma 2}~\cite{BGP}.
Given an algebra $A$, an RB-operator $P$ on $A$ of weight~$\lambda$,
and an (anti)automorphism of~$A$,
the operator $P^{(\psi)} = \psi^{-1}P\psi$
is again an RB-operator of weight~$\lambda$ on~$A$.

{\bf Lemma 3}~\cite{BGP,GuoMonograph}.
Let $A$ be a unital algebra, $R$ be an RB-operator of weight~0 on~$A$.

(a) $1\not\in \Imm(R)$;

(b) if $R(1) = 0$, then $\Imm(R)\subseteq \ker(R)$ and $R^2 = 0$;

(c) $(R(1))^n = n! R^n(1)$, $n\in\mathbb{N}$.

Now, we give not so simple results concerning RB-operators defined on unital algebras.

{\bf Lemma 4}~\cite{Unital}.
Let $A$ be a unital associative finite-dimensional algebra and $R$ be an RB-operator
on $A$ of weight zero. Then $R(1)$ is nilpotent.

{\bf Theorem 1}~\cite{Unital,Miller66}.
Let $A$ be a unital associative finite-dimensional algebra over a~field of characteristic zero.
Then there exists $N$ such that $R^N = 0$
for every RB-operator~$R$ of weight~0 on $A$.

The main object for the study is the algebra~$U_3(F)$ of upper-triangular matrices of order~3.
We write down an arbitrary automorphism $\psi$ of $U_3(F)$~\cite{Gub2021}:
\begin{equation}\label{AutU}
\begin{gathered}
e_{11}\to \begin{pmatrix}
1 & \beta & \gamma \\
0 & 0 & 0 \\
0 & 0 & 0 \\
\end{pmatrix},\quad
e_{12}\to \begin{pmatrix}
0 & \delta & \varepsilon \\
0 & 0 & 0 \\
0 & 0 & 0 \\
\end{pmatrix},\quad
e_{13}\to \begin{pmatrix}
0 & 0 & \alpha \\
0 & 0 & 0 \\
0 & 0 & 0 \\
\end{pmatrix}, \\
e_{22}\to
\begin{pmatrix}
0 & -\beta & -\beta\varepsilon/\delta \\
0 & 1 & \varepsilon/\delta \\
0 & 0 & 0 \\
\end{pmatrix}\!, \quad
e_{23}\to
\frac{1}{\delta}\begin{pmatrix}
0 & 0 & -\alpha\beta \\
0 & 0 & \alpha \\
0 & 0 & 0 \\
\end{pmatrix}\!,\quad
e_{33}\to
\begin{pmatrix}
0 & 0 & \beta\varepsilon/\delta -\gamma \\
0 & 0 & -\varepsilon/\delta \\
0 & 0 & 1 \\
\end{pmatrix}\!,
\end{gathered}
\end{equation}
where $\alpha,\delta \neq0$.

To apply Lemma~2, we also give an antiautomorphism~$\Theta_{1,3}$ of~$U_3(F)$ defined as follows,
\begin{equation}\label{Theta13}
e_{11} \to e_{33}, \quad
e_{12} \to e_{23}, \quad
e_{13} \to e_{13}, \quad
e_{22} \to e_{22}, \quad
e_{23} \to e_{12}, \quad
e_{33} \to e_{11}.
\end{equation}
Equivalently, $\Theta_{1,3}$ is a restriction on~$U_3(F)$ of the antiautomorphism of~$M_3(F)$ defined by the formula
$\Theta_{1,3}(X) = ZX^T Z$, where
$Z = Z^{-1} = \begin{pmatrix}
0 & 0 & 1 \\
0 & 1 & 0 \\
1 & 0 & 0 \\
\end{pmatrix}$ and $T$ denotes the transpose of a matrix.

Let us give a sufficiently general example of RB-operator of weight~0.

{\bf Example}~\cite{Unital}.
Let $A$ be an algebra such that
$A = B\oplus C$ (as vector spaces) with abelian subalgebra $C$,
moreover, $C$ is a $B$-bimodule.
Then any linear operator~$R$ on~$A$ defined as follows,
$R\colon B\to C$, $R\colon C\to (0)$ is an RB-operator of weight~0 on~$A$.

For $A = U_3(F)$, Example provides that the following operators are indeed Rota---Baxter ones:
\begin{gather*}
R(e_{12}) = R(e_{13}) = 0, \quad
R(e_{11}),
R(e_{22}),
R(e_{33}),
R(e_{23})\in L(e_{12},e_{13}); \\
R(e_{13}) = 0, \quad
R(e_{11}), R(e_{12}), R(e_{22}), R(e_{23}), R(e_{33})\in L(e_{13}).
\end{gather*}

\section{Strategy of the solution}

Let $R$ be an RB-operator on~$U_3(F)$ and $\charr F = 0$.
By~Lemma~4, $R(1)$ is nilpotent.
Due to~\cite{Belitskii}, we know that up to the action of antiautomorphism of~$U_3(F)$, we have
the following four variants:
\begin{itemize}
 \item $R(1) = 0$,
 \item $R(1) = e_{12}$,
 \item $R(1) = e_{13}$,
 \item $R(1) = e_{12}+e_{23}$.
\end{itemize}
By Lemma~3c), $2R^2(1) = R(1)^2$, so we additionally derive that
$R^2(1) = 0$ in the second and third cases, and $R(1)^2 = e_{13}/2$ in the fourth one. 
Such knowledge is sufficient to start the computations and then deal with all possible subcases.

When $R(1) = 0$, we still need 30 parameters to describe an RB-operator. However, an important conclusion helps: $R^2 = 0$ (Lemma~3b).
Thus, we separate this most difficult case into subcases.
We study the semisimple part of the Malcev---Wedderburn decomposition of the associative algebra~$\Imm(R)$.

If $\Imm(R)$ is nilpotent, then $\Imm(R) \subseteq V := L(e_{12},e_{13},e_{23})$, and it is not hard to finish this variant involving only 15 parameters.

If $\Imm(R)\cong F\dotplus U$, where $U$ is the (nilpotent) radical of~$\Imm(R)$, then we clarify, what a matrix generates the semisimple part of~$\Imm(R)$.

{\bf Lemma 5}.
Let $A$ be an idempotent matrix of~$U_3(F)$.

a) If $\rank(A) = 1$, then up to the action of an antiautomorphism of~$U_3(F)$,
we have either $A = e_{11}$ or $A = e_{22}$.

b) If $\rank(A) = 2$, then up to the action an antiautomorphism of~$U_3(F)$,
we have either $A = e_{11}+e_{22}$ or $A = e_{11}+e_{33}$.

{\sc Proof}.
a) Note that the conditions of the statement imply that only one of the diagonal elements~$a_{ii}$ of $A$ is nonzero, which means that $a_{ii} = 1$.
Suppose that $a_{11} = 1$. Since $A^2 = A$, the matrix $A$ has the form $e_{11}+be_{12}+ce_{13}$. It remains to apply~$\psi^{-1}$ (see~\eqref{AutU}) with $\beta = b$ and $\gamma = c$ to get~$e_{11}$.
If $a_{22} = 1$, we analogously apply~$\psi$ and get~$e_{22}$.
If $a_{33} = 1$, we apply~$\Theta_{1,3}$ (see~\eqref{Theta13}) and come to the already studied case~$a_{11} = 1$.

b) Suppose that $a_{11} = a_{22} = 1$, other cases may be considered similarly.
Thus, $A = e_{11} + be_{12} + ce_{13}+ e_{22} + d e_{23}$.
Now we use~$\psi^{-1}$ with $(\alpha,\beta,\gamma,\delta,\varepsilon) = (1,b,c,1,d)$ and get~$e_{11}+e_{22}$.
\hfill $\square$

Thus, we separately consider all four subcases

\noindent--- $\Imm(R) = L(e_{11}) + U$,

\noindent--- $\Imm(R) = L(e_{22}) + U$,

\noindent--- $\Imm(R) = L(e_{11}+e_{22}) + U$,

\noindent--- $\Imm(R) = L(e_{11}+e_{33}) + U$,
where $U\subseteq V$.

Since $R^2 = 0$, we conclude that $R(e_{11}) = 0$, $R(e_{22}) = 0$, $R(e_{11}+e_{22}) = 0$ or $R(e_{11}+e_{33}) = 0$ respectively. 
In each of the subcases, we need 16 parameters, and the computations allow us to deal with them.

Finally, we may have $\Imm(R)\cong (F\oplus F)\dotplus U$ for some $U\subseteq V$. Denote by $A$ and~$B$ idempotents of~$\Imm(R)$ such that $AB = 0$ and $A+B$ is the unit of the semisimple part of~$\Imm(R)$.
Since $A+B$ is idempotent, we have $a_{ii}b_{ii} = 0$ for $i=1,2,3$ (otherwise we would obtain $a_{ii}+b_{ii} = 2$).
Applying the property that $\Imm(R)$ does not contain a non-degenerate matrix (it follows from Lemma 3a), we get that $\rank(A) = \rank(B) = 1$.
By Lemma~5, we may assume that~$A = e_{11}$ or $A = e_{22}$.
Let us consider the first variant; the second one is analogous.
If $b_{22} = 1$, then $AB = 0$ implies that $B = e_{22} + k e_{23}$. Now, $\psi^{-1}\in\Aut(F)$ with $\varepsilon = k$, $\alpha = \delta = 1$, $\beta = \gamma = 0$, maps $\Imm(R)$ to $L(e_{11},e_{22})\oplus U'$, $U'\subseteq V$.
If $b_{33} = 1$, analogously, we may assume that $\Imm(R)=L(e_{11},e_{33})\oplus U'$, $U'\subseteq V$.
If $A = e_{22}$, then with the help of (anti)automorphisms of~$U_3(F)$, we come to the variant~$\Imm(R) = L(e_{11},e_{22})\oplus U'$, $U'\subseteq V$.

Let us highlight that below we fulfill computations even with some slightly mild conditions that we have mentioned above.
For example, when $R(1) = 0$, $\Imm(R) = L(e_{11}) + U$ for $U\subseteq V$, we apply that $R(e_{11}) = 0$. However,
we do not try to put the fact that there exists~$A\in U_3(F)$ such that $R(A) = e_{11}$ in the computational process.
We find it simpler to avoid the search for such~$A$ and focus on the restrictions of the parameters~$b_{ij,kl}$ that arise from
$$
R(e_{ij}) = \sum\limits_{k,l=1,\,k\leq l}^3 b_{ij,kl}e_{kl},\quad e_{ij}\in U_3(F).
$$
Surprisingly, such an approach is enough to solve the general problem.

\section{Case $R(1) = 0$}

By Lemma~3b), the condition $R(1)=0$ implies $R^2 = 0$.
We split the general case $R(1) = 0$ on the subcases depending on which diagonal matrices lie in $\Imm(R)$.

\subsection{$\Imm(R)$ is nilpotent}

Computations with \texttt{Singular} imply $R(e_{13}) = R(1) = 0$ and
$$
R(e_{12}) = \begin{pmatrix}
0 & 0 & a \\
0 & 0 & b \\
0 & 0 & 0
\end{pmatrix}\!, \
R(e_{23}) = \begin{pmatrix}
0 & c & d \\
0 & 0 & 0 \\
0 & 0 & 0
\end{pmatrix}\!, \
R(e_{22}) = \begin{pmatrix}
0 & e & f \\
0 & 0 & g \\
0 & 0 & 0
\end{pmatrix}\!,\
R(e_{33}) = \begin{pmatrix}
0 & h & i \\
0 & 0 & j \\
0 & 0 & 0
\end{pmatrix}\!\!
$$
joint with the quadratic relations
$$
\{ c, e, h\} \times a = 0, \quad
\{ c, d, e, h \} \times \{ b, g, j \} = 0.
$$

{\sc Subcase 1}: $(c,e,h)\neq(0,0,0)$.
Then $a = b = g = j = 0$, and we get an RB-operator
$$
R(e_{12}) = R(e_{13}) = R(1) = 0, \quad
R(e_{22}), R(e_{23}), R(e_{33})\in L(e_{12},e_{13}).
$$
Here we may omit the condition that not all of $c,e,h$ are nonzero, since $R$ is an RB-operator in any case.

{\sc Subcase 2}: $c = e = h = 0$. It remains to consider the last fork.

{\sc Subcase 2.1}: $d\neq0$.
Then $b = g = j = 0$, it is an RB-operator
$$
R(e_{13}) = R(1) = 0, \quad
R(e_{12}), R(e_{22}), R(e_{23}), R(e_{33}) \in L(e_{13}).
$$

{\sc Subcase 2.2}: $d = 0$.
Hence, $R$ has the form
$$
R(e_{13}) = R(e_{23}) = R(1) = 0, \quad
R(e_{12}), R(e_{22}), R(e_{33})\in L(e_{13},e_{23}),
$$
which, up to conjugation with~$\Theta_{1,3}$, coincides with the RB-operator from Subcase~1.

\subsection{Semisimple part of $\Imm(R)$ is $L(e_{11})$}

Since $R^2 = 0$, we have $R(e_{11}) = 0$.
Computations with \texttt{Singular} imply
\begin{gather*}
R(e_{12}) = \begin{pmatrix}
e & -a & f \\
0 & 0 & g \\
0 & 0 & 0
\end{pmatrix}, \quad
R(e_{13}) = \begin{pmatrix}
k & h & a \\
0 & 0 & l \\
0 & 0 & 0
\end{pmatrix}, \quad
R(e_{23}) = \begin{pmatrix}
h & i & j \\
0 & 0 & 0 \\
0 & 0 & 0
\end{pmatrix}, \\
R(e_{22}) = \begin{pmatrix}
-a & -b & -c \\
0 & 0 & -d \\
0 & 0 & 0
\end{pmatrix}, \quad
R(e_{33}) = \begin{pmatrix}
a & b & c \\
0 & 0 & d \\
0 & 0 & 0
\end{pmatrix}
\end{gather*}
joint with the quadratic relations, including
$$
\{ d, g, l \} \times \{ a, b, h, i, j, k \} = 0.
$$

{\sc Subcase 1}: $(d,g,l)\neq(0,0,0)$.
Then $a = b = h = i = j = k = 0$,
and the following relations hold:
$$
\{ c, d, f, l+e \} \times l = 0, \quad
ce = 0, \quad
g(l+e) = 0.
$$

{\sc Subcase 1.1}: $l\neq0$.
Hence, $c = d = f = 0$ and $l = -e\neq0$.
Dividing by~$l$ and applying~$\psi$ with $\alpha = \delta = 1$, $\beta = \gamma = 0$ and $\varepsilon = -g/l$, we get an RB-operator
\begin{equation}
R(e_{11}) = R(e_{22}) = R(e_{23}) = R(e_{33}) = 0, \quad
R(e_{12}) = - e_{11}, \quad
R(e_{13}) = e_{23}.
\end{equation}

{\sc Subcase 1.2}: $l = 0$.
It means that $ce = ge = 0$.
Since the projection of $\Imm(R)$ on $e_{11}$ is nonzero, we have $e\neq0$.
Thus, we get $c = g = 0$ and an RB-operator
$$
R(e_{11}) = R(e_{13}) = R(e_{23}) = R(1) = 0, \quad
R(e_{12}) = e e_{11} + f e_{13}, \quad
R(e_{33}) = de_{23}, \quad e,d\neq0.
$$
An automorphism~$\psi$ of $U_3(F)$ defined with
$\alpha = d/e$, $\delta = 1/e$, $\gamma = f/e$, $\beta = \varepsilon = 0$ gives an RB-operator
\begin{equation}
R(e_{11}) = R(e_{13}) = R(e_{23}) = 0, \quad
R(e_{12}) =  e_{11}, \quad
R(e_{22}) = -e_{23}, \quad
R(e_{33}) =  e_{23}.
\end{equation}

{\sc Subcase 2}: $d = g = l = 0$.
Then
$R(e_{11}) = 0$ and
\begin{gather*}
R(e_{12}) = \begin{pmatrix}
e & -a & f \\
0 & 0 & 0 \\
0 & 0 & 0
\end{pmatrix}, \quad
R(e_{13}) = \begin{pmatrix}
k & h & a \\
0 & 0 & 0 \\
0 & 0 & 0
\end{pmatrix}, \quad
R(e_{23}) = \begin{pmatrix}
h & i & j \\
0 & 0 & 0 \\
0 & 0 & 0
\end{pmatrix}, \\
R(e_{22}) = \begin{pmatrix}
-a & -b & -c \\
0 & 0 & 0 \\
0 & 0 & 0
\end{pmatrix}, \quad
R(e_{33}) = \begin{pmatrix}
a & b & c \\
0 & 0 & 0 \\
0 & 0 & 0
\end{pmatrix}.
\end{gather*}
Moreover, the coefficients satisfy the following quadratic relations:
\begin{gather*}
\{ a, e, h,k \} \times (j-b) = 0, \quad
ab = ch, \quad
fb = -ac, \quad
eb = -kc, \quad
a^2 = kc, \\
ia = hb, \quad
ha = kb, \quad
fa = ec, \quad
fi = -hc, \\
ei = -kb, \quad
h^2 = ki, \quad
fh = -kc, \quad
eh = -ka, \quad
fk = ea.
\end{gather*}

{\sc Subcase 2.1}: $k = 0$.
Then $a = h = 0$.
Since $e_{11}\in\Imm(R)$, we deduce that $e\neq0$
and $b = c = i = j = 0$.
Dividing by $e$ and applying a corresponding automorphism, we get an RB-operator
\begin{equation}
R(e_{11}) = R(e_{22}) = R(e_{33}) = R(e_{13}) =  R(e_{23}) = 0, \quad
R(e_{12}) = e_{11}.
\end{equation}

{\sc Subcase 2.2}: $k \neq 0$.
Then $b = j$ and we may express
$$
a = -eh/k,\quad
b = -eh^2/k^2, \quad
c = e^2h^2/k^3, \quad
i = h^2/k, \quad
f = -e^2h/k^2.
$$
Hence, after dividing $R$ by~$k$, we get an RB-operator
$R(e_{11}) = R(1) = 0$ and for $s = e/k$ and $t = h/k$
\begin{gather*}
R(e_{12}) = \begin{pmatrix}
s & st & -s^2t \\
0 & 0 & 0 \\
0 & 0 & 0
\end{pmatrix}, \quad
R(e_{13}) = \begin{pmatrix}
1 & t & -st \\
0 & 0 & 0 \\
0 & 0 & 0
\end{pmatrix}, \\
R(e_{23}) = \begin{pmatrix}
t & t^2 & -st^2 \\
0 & 0 & 0 \\
0 & 0 & 0
\end{pmatrix}, \quad
R(e_{33}) = \begin{pmatrix}
-st & -st^2 & s^2t^2 \\
0 & 0 & 0 \\
0 & 0 & 0
\end{pmatrix}.
\end{gather*}
Applying~$\psi$ with $\alpha = \delta = 1$, $\beta = t$, $\gamma = -st$, $\varepsilon = -s$, we get an RB-operator
\begin{equation}
R(e_{11}) = R(e_{22}) = R(e_{33}) = R(e_{12}) = R(e_{23}) = 0, \quad
R(e_{13}) =  e_{11}.
\end{equation}

\subsection{Semisimple part of $\Imm(R)$ is $L(e_{22})$}

We have $R(e_{22}) = 0$.
Moreover, computations with \texttt{Singular} imply \begin{gather*}
R(e_{12}) = \begin{pmatrix}
0 & 0 & d \\
0 & e & f \\
0 & 0 & 0
\end{pmatrix}, \quad
R(e_{13}) = \begin{pmatrix}
0 & k & 0 \\
0 & 0 & l \\
0 & 0 & 0
\end{pmatrix}, \quad
R(e_{23}) = \begin{pmatrix}
0 & i & j \\
0 & h & 0 \\
0 & 0 & 0
\end{pmatrix}, \\
R(e_{11}) = \begin{pmatrix}
0 & -a & -b \\
0 & 0 & -c \\
0 & 0 & 0
\end{pmatrix}, \quad
R(e_{33}) = \begin{pmatrix}
0 & a & b \\
0 & 0 & c \\
0 & 0 & 0
\end{pmatrix}.
\end{gather*}
joint with the quadratic relations, including
$$
\{ a, h, i, k \} \times \{ c, d, e, f, l \} = 0.
$$

Since $e_{22}\in\Imm(R)$, one of $e,h$ is nonzero.
Up to conjugation with $\Theta_{13}$, we may assume that $h\neq0$.
Hence, $c = d = e = f = l$, and the rest of the quadratic equations are
$$
\{ b, i, j, h-k \} \times k = 0, \quad
hj = 0, \quad
a(k-h) = 0.
$$
We have $h\neq0$, hence, $j = 0$.

{\sc Subcase 1}: $k\neq0$.
Then $b = i = 0$ and $h = k$.
After dividing by~$k$, we get an RB-operator
$$
R(e_{12}) = R(e_{22}) = 0, \quad
R(e_{11}) = -se_{12}, \quad
R(e_{33}) = se_{12}, \quad
R(e_{13}) = e_{12}, \quad
R(e_{23}) = e_{22}.
$$
Applying~$\psi\in\Aut(U_3(F))$ with $\alpha = \delta = 1$, $\beta = \varepsilon = 0$, $\gamma = s$, we get~$R$ of the form
\begin{equation}
R(e_{12}) = R(e_{11}) = R(e_{22}) = R(e_{33}) = 0, \quad
R(e_{13}) = e_{12}, \quad
R(e_{23}) = e_{22}.
\end{equation}

{\sc Subcase 2}: $k = 0$.
Since $e_{22}\in\Imm(R)$, $h\neq0$, it implies that $a = 0$, and $R$ equals after dividing by~$h$:
$$
R(e_{12}) = R(e_{13}) = R(e_{22}) = 0, \quad
R(e_{11}) = -se_{13}, \quad
R(e_{33}) = se_{13}, \quad
R(e_{23}) = e_{22} + te_{12}.
$$
If $s = 0$, then conjugation with~$\psi$ with $\alpha = \delta = 1$, $\beta = -t$, $\gamma = \varepsilon = 0$ implies
\begin{equation}
R(e_{12}) = R(e_{13}) = R(e_{11}) = R(e_{22}) = R(e_{33}) = 0, \quad
R(e_{23}) = e_{22}.
\end{equation}
If $s\neq0$, then applying~$\psi$ with $\alpha = \delta = s$, $\beta = -t$, $\gamma = \varepsilon = 0$, we get
\begin{equation}
R(e_{12}) = R(e_{13}) = R(e_{22}) = 0, \quad
R(e_{11}) = -e_{13}, \quad
R(e_{33}) = e_{13},\quad
R(e_{23}) = e_{22}.
\end{equation}

\subsection{Semisimple part of $\Imm(R)$ is $L(e_{11}+e_{22})$}

We have $R(e_{33}) = 0$.
Computations with \texttt{Singular} imply
\begin{gather*}
R(e_{12}) = \begin{pmatrix}
-h & 0 & k \\
0 & -h & l \\
0 & 0 & 0
\end{pmatrix}, \quad
R(e_{13}) = \begin{pmatrix}
i & j & -g \\
0 & i & h \\
0 & 0 & 0
\end{pmatrix}, \quad
R(e_{23}) = \begin{pmatrix}
d & e & f \\
0 & d & g \\
0 & 0 & 0
\end{pmatrix}, \\
R(e_{11}) = \begin{pmatrix}
g & -a & -b \\
0 & g & -c \\
0 & 0 & 0
\end{pmatrix}, \quad
R(e_{22}) = \begin{pmatrix}
-g & a & b \\
0 & -g & c \\
0 & 0 & 0
\end{pmatrix}.
\end{gather*}
joint with the quadratic relations, including
$$
\{ a,e,j \} \times \{ c,g,h,i,k,l\} = 0.
$$

{\sc Subcase 1}: $(a,e,j)\neq(0,0,0)$.
Then $c = g = h = i = k = l = 0$, and the complete list of quadratic equations is
$$
\{ a,b,f,j-d \} \times j = 0, \quad
e(j-d) = 0, \quad
bd = 0.
$$
Since $e_{11}+e_{22}\in\Imm(R)$, we have $d\neq0$.
Thus, $b = 0$.

{\sc Subcase 1.1}: $j\neq d$.
Then $e = j = 0$. After dividing by~$d$ we get
$$
R(e_{12}) = R(e_{13}) = R(e_{33}) = 0, \
R(e_{23}) = e_{11}+e_{22} + te_{13},\
R(e_{11}) = se_{12},\
R(e_{22}) = -se_{12}.
$$

Since $e_{11}+e_{22}\in \Imm(R)$, we derive $t = 0$.
We apply conjugation~$\psi$ with $\alpha = \delta \neq 0$, $\beta = \gamma =  \varepsilon = 0$, to get, depending on if $s = 0$ or $s\neq0$, the following two operators:
\begin{gather}
R(e_{12}) = R(e_{13}) = R(e_{11}) = R(e_{22}) = R(e_{33}) = 0, \quad
R(e_{23}) = e_{11}+e_{22}; \label{11+22-first} \\
R(e_{12}) = R(e_{13}) = R(e_{33}) = 0, \
R(e_{23}) = e_{11}+e_{22},\
R(e_{11}) = e_{12},\
R(e_{22}) = -e_{12},
\end{gather}
where we take $\alpha = \delta = s$ in the last one.

{\sc Subcase 1.2}: $j = d \neq0$. Hence, $a = f = 0$.
Dividing by~$d$, we come to
$$
R(1) = R(e_{12}) = R(e_{33}) = 0, \quad
R(e_{13}) = e_{12}, \quad
R(e_{23}) = e_{11} + e_{22} + se_{12}, \quad
R(e_{22}) = te_{13}.
$$
By the conditions, $s = 0$.
Depending on the value~$t$ and applying for $t\neq0$ the corresponding~$\psi$, we have two RB-operators:
\begin{gather}
R(e_{12}) = R(e_{11}) = R(e_{22}) = R(e_{33}) = 0, \quad
R(e_{13}) = e_{12}, \quad
R(e_{23}) = e_{11} + e_{22};\\
R(e_{12}) {=} R(e_{33}) {=} 0, \
R(e_{13}) = e_{12}, \
R(e_{23}) = e_{11} {+} e_{22}, \
R(e_{11}) = e_{13}, \
R(e_{22}) = -e_{13}.
\end{gather}

{\sc Subcase 2}: $a = e = j = 0$.
Among the relations, we have the following ones: $\{d,i,g,h\} \times (c+k) = 0$.
Since $e_{11}+e_{22}\in\Imm(R)$, we conclude that $k = -c$.
Computations provide a system of quadratic equations:
\begin{gather*}
g^2 + cd = 0, \quad
bd + gf = 0 \quad
hd = ig, \quad
ld + hg = 0, \quad
cf = bg, \\
hf + g^2 = 0, \quad
if + dg = 0, \quad
lf + cg = 0, \quad
hb + cg = 0, \quad
ib = g^2, \\
lg - ch = 0, \quad
li + h^2 = 0, \quad
ci + hg = 0.
\end{gather*}

If $i = 0$, then $g = h = 0$.
Since $d\neq0$ (otherwise $e_{11}+e_{22}\not\in\Imm(R)$),
we also have $b = c = l = 0$, it is already obtained RB-operator~\eqref{11+22-first}.

Consider $i\neq0$. Therefore, we express
$$
b = d^2h^2/i^3, \quad
c = -dh^2/i^2, \quad
f = -d^2h/i^2, \quad
g = dh/i, \quad
l = -h^2/i.
$$
Dividing by~$i$, we get in terms of $s = d/i$ and $t = h/i$
that $R(1) = R(e_{33}) = 0$,
\begin{gather*}
R(e_{12}) = \begin{pmatrix}
-s & 0 & s^2t \\
0 & -s & -s^2 \\
0 & 0 & 0
\end{pmatrix}, \quad
R(e_{13}) = \begin{pmatrix}
1 & 0 & -st \\
0 & 1 & s \\
0 & 0 & 0
\end{pmatrix}, \\
R(e_{23}) = \begin{pmatrix}
t & 0 & -st^2 \\
0 & t & st \\
0 & 0 & 0
\end{pmatrix}, \quad
R(e_{22}) = \begin{pmatrix}
-st & 0 & s^2t^2 \\
0 & -st & -s^2t \\
0 & 0 & 0
\end{pmatrix}.
\end{gather*}
Applying conjugation with~$\psi$ defined with parameters
$\alpha = \delta = 1$, $\beta = t$, $\gamma = 0$, $\varepsilon = s$, we get
\begin{equation}
R(e_{11}) = R(e_{22}) = R(e_{12}) = R(e_{23}) = R(e_{33}) = 0, \quad
R(e_{13}) = e_{11}+e_{22}.
\end{equation}

\subsection{Semisimple part of $\Imm(R)$ is $L(e_{11}+e_{33})$}

We have $R(1) = R(e_{22}) = 0$.
Computations with \texttt{Singular} imply
\begin{gather*}
R(e_{12}) = \begin{pmatrix}
k & 0 & l \\
0 & 0 & m \\
0 & 0 & k
\end{pmatrix}, \quad
R(e_{13}) = \begin{pmatrix}
0 & g & 0 \\
0 & 0 & h \\
0 & 0 & 0
\end{pmatrix}, \quad
R(e_{23}) = \begin{pmatrix}
d & e & f \\
0 & 0 & 0 \\
0 & 0 & d
\end{pmatrix}, \\
R(e_{11}) = \begin{pmatrix}
0 & -a & -b \\
0 & 0 & -c \\
0 & 0 & 0
\end{pmatrix}, \quad
R(e_{33}) = \begin{pmatrix}
0 & a & b \\
0 & 0 & c \\
0 & 0 & 0
\end{pmatrix}.
\end{gather*}
joint with the quadratic relations, including
$\{d,k\} \times \{f, l\} = 0$.
Since $e_{11}+e_{33}\in\Imm(R)$, we have $(d,k)\neq(0,0)$.
Thus, $f = l = 0$.
Among the new list of quadratic relations, we find $dk = 0$.
Up to conjugation with $\Theta_{1,3}$, we may assume that $d = 0$. It implies $g = 0$ and the following relations:
$$
\{ b,e,m,k+h \} \times h = 0, \quad
\{ a,e,k+h\} \times c = 0, \quad
\{ h, k, m \} \times a = 0, \quad
\{ k, m \} \times e = 0.
$$
Since $k \neq 0$, we conclude that $a = e = 0$.
The rest of the equations are
$$
\{ b,m,k+h \} \times h = 0, \quad
c(k+h) = 0.
$$

{\sc Subcase 1}: $k+h\neq 0$.
Then $c = h = 0$ and after dividing by~$k$ we get
$$
R(1) = R(e_{22}) = R(e_{13}) = R(e_{23}) = 0, \quad
R(e_{12}) = e_{11} + e_{33} + se_{23}, \quad
R(e_{33}) = te_{13}.
$$
Conjugation with $\psi$ defined with $\beta = \gamma = 0$, $\delta = 1$, $\varepsilon = -s$, $\alpha\neq0$ gives, 
depending on whether $t = 0$ or $t\neq0$, one of the following variants:
\begin{gather}
R(e_{11}) = R(e_{22}) = R(e_{13}) = R(e_{23}) = R(e_{33}) = 0, \quad
R(e_{12}) = e_{11} + e_{33}; \\
R(e_{22}) = R(e_{13}) = R(e_{23}) = 0, \quad
R(e_{12}) = e_{11} + e_{33}, \
R(e_{11}) = -e_{13}, \
R(e_{33}) = e_{13},
\end{gather}
in the second one, we take $\alpha = t$.

{\sc Subcase 2}: $k+h = 0$. Hence, $h = -k\neq0$
and so, $b = m = 0$.
Dividing by~$k$, write down
$$
R(1) = R(e_{22}) = R(e_{23}) = 0, \quad
R(e_{12}) = e_{11} + e_{33} + se_{23}, \quad
R(e_{13}) = -e_{23}, \quad
R(e_{33}) = te_{23}.
$$
Again, conjugation with $\psi$ defined with $\beta = \gamma = 0$, $\delta = 1$, $\varepsilon = -s$, $\alpha\neq0$ gives, 
depending on whether $t = 0$ or $t\neq0$, one of the following variants:
\begin{gather}
R(e_{11}) = R(e_{22}) = R(e_{23}) = R(e_{33}) = 0, \quad
R(e_{13}) = -e_{23}, \quad
R(e_{12}) = e_{11} + e_{33}; \\
R(e_{22}) {=} R(e_{23}) {=} 0, \
R(e_{13}) {=} -e_{23}, \
R(e_{12}) {=} e_{11} {+} e_{33}, \
R(e_{11}) = -e_{23}, \
R(e_{33}) = e_{23},
\end{gather}
in the second one, we take $\alpha = t$ too.

\subsection{Semisimple part of $\Imm(R)$ is $L(e_{11},e_{22})$}

Computations with \texttt{Singular} imply
$R(e_{11}) = R(e_{22}) = R(e_{33}) = 0$ and
$$
R(e_{12}) = \begin{pmatrix}
i & 0 & j \\
0 & k & l \\
0 & 0 & 0
\end{pmatrix}, \quad
R(e_{13}) = \begin{pmatrix}
e & f & 0 \\
0 & g & h \\
0 & 0 & 0
\end{pmatrix}, \quad
R(e_{23}) = \begin{pmatrix}
a & b & c \\
0 & d & 0 \\
0 & 0 & 0
\end{pmatrix},
$$
satisfying additional quadratic relations, including
$$
\{ a,b,c,d\} \times \{h,k,l\} = 0.
$$

{\sc Subcase 1}: $(h,k,l)\neq(0,0,0)$.
Then $R(e_{11}) = R(e_{22}) = R(e_{33}) = R(e_{23}) = 0$, and
the computations show that $f = 0$, the coefficients of $R$ satisfy
the quadratic relations, including the following ones:
$$
\{ g, h, k \} \times \{ j, e-g \} = 0.
$$

If $(j,e-g)\neq(0,0)$, then $g = h = k = 0$ and so, $e_{22}\not\in \Imm(R)$.
Hence, $j = 0$ and $e = g$, and the following quadratic equations hold:
$$
\{ g, l, h-k \} \times (i+h) = 0, \quad
lg + ik = 0.
$$
The restrictions on $\Imm(R)$ imply that $g\neq0$, so,
$i = -h$ and $l = kh/g$.
After dividing by~$g$, we get for $s = k/g$ and $t = h/g$
$$
R(e_{12}) = \begin{pmatrix}
-t & 0 & 0 \\
0 & s & st \\
0 & 0 & 0
\end{pmatrix}, \quad
R(e_{13}) = \begin{pmatrix}
1 & 0 & 0 \\
0 & 1 & t \\
0 & 0 & 0
\end{pmatrix}, \quad s+t\neq0.
$$
Applying~$\psi$ with $\alpha = 1$, $\beta = \gamma = 0$,
$\delta = 1/(s+t)$ and $\varepsilon = t/(s+t)$, we get an RB-operator
\begin{equation}
R(e_{11}) = R(e_{22}) = R(e_{33}) = R(e_{23}) = 0, \quad
R(e_{12}) = e_{22}, \quad
R(e_{13}) = e_{11} + e_{22}.
\end{equation}

{\sc Subcase 2}: $h = k = l = 0$.
Computations with \texttt{Singular} imply the following relations:
$$
\{ i,j \} \times \{ a-d, b, f, g \} = 0.
$$

{\sc Subcase 2.1}: $(i,j)\neq(0,0)$.
Hence, $b = f = g = 0$ and $a = d$. We get
$R(e_{11}) = R(e_{22}) = R(e_{33}) = 0$ and
$$
R(e_{12}) = \begin{pmatrix}
i & 0 & j \\
0 & 0 & 0 \\
0 & 0 & 0
\end{pmatrix}, \quad
R(e_{13}) = \begin{pmatrix}
e & 0 & 0 \\
0 & 0 & 0 \\
0 & 0 & 0
\end{pmatrix}, \quad
R(e_{23}) = \begin{pmatrix}
a & 0 & c \\
0 & a & 0 \\
0 & 0 & 0
\end{pmatrix},
$$
satisfying the four quadratic equations
$$
\{ a, c, j \} \times e = 0, \quad
ic = ja.
$$
Since $\Imm(R)$ contains $e_{11},e_{22}$, we conclude that $a\neq0$. Therefore, $e = 0$, $j = ic/a$, and $a,i\neq0$.
Dividing~$R$ by~$a$, we get for $s = i/a$ and $t = c/a$
$$
R(e_{12}) = s(e_{11}+te_{13}),\quad
R(e_{23}) = e_{11}+e_{22}+te_{13},\quad s\neq0.
$$
Applying~$\psi$ with $\alpha = \delta = 1/s$, $\beta = \varepsilon = 0$, $\gamma = t$, we get an RB-operator
\begin{equation}
R(e_{11}) = R(e_{22}) = R(e_{33}) = R(e_{13}) = 0, \quad
R(e_{12}) = e_{11}, \quad
R(e_{23}) = e_{11} + e_{22}.
\end{equation}

{\sc Subcase 2.2}: $i = j = 0$.
Thus, we have
$R(e_{11}) = R(e_{22}) = R(e_{33}) = R(e_{12}) = 0$ and
$$
R(e_{13}) = \begin{pmatrix}
e & f & 0 \\
0 & g & 0 \\
0 & 0 & 0
\end{pmatrix}, \quad
R(e_{23}) = \begin{pmatrix}
a & b & c \\
0 & d & 0 \\
0 & 0 & 0
\end{pmatrix},
$$
satisfying due to \texttt{Singular}
the quadratic relations, including
$$
\{ e-g, a-d\} \times g = 0, \quad ce = 0.
$$
If $g\neq0$, then $e = g$ and $a = d$, a contradiction to the conditions on $\Imm(R)$. Thus, $g = 0$ and $d,e\neq0$.
By this reason, the equation $ce = 0$ implies that $c = 0$. It remains to consider the equations
$$
\{ b, e, a+f-d \} \times (a-f) = 0, \quad
a^2 - ad - be = 0.
$$
Again, since $e\neq0$, we conclude that $a = f$ and $b = a(a-d)/e$.
After diving on~$e$, we get for $s = a/e$ and $t = d/e\neq0$
$$
R(e_{13}) = e_{11} + se_{12}, \quad
R(e_{23}) = s( e_{11} + (s-t)e_{12}) + te_{22}.
$$
Applying~$\psi$ with $\alpha = 1$, $\delta = t$, $\gamma = \varepsilon = 0$, $\beta = s$, we get an RB-operator
\begin{equation}
R(e_{11}) = R(e_{22}) = R(e_{33}) = R(e_{12}) = 0, \quad
R(e_{13}) = e_{11}, \quad
R(e_{23}) = e_{22}.
\end{equation}

\subsection{Semisimple part of $\Imm(R)$ is $L(e_{11},e_{33})$}

Computations with \texttt{Singular} imply
$R(e_{11}) = R(e_{22}) = R(e_{33}) = 0$ and
$$
R(e_{12}) = \begin{pmatrix}
i & 0 & j \\
0 & 0 & l \\
0 & 0 & k
\end{pmatrix}, \quad
R(e_{13}) = \begin{pmatrix}
e & f & 0 \\
0 & 0 & h \\
0 & 0 & g
\end{pmatrix}, \quad
R(e_{23}) = \begin{pmatrix}
a & b & c \\
0 & 0 & 0 \\
0 & 0 & d
\end{pmatrix},
$$
satisfying additional quadratic relations, including
$$
\{ a,b,c,d\} \times \{h,k,l\} = 0.
$$

{\sc Subcase 1}: $(h,k,l)\neq(0,0,0)$.
Then $R(e_{11}) = R(e_{22}) = R(e_{33}) = R(e_{23}) = 0$; moreover,
\texttt{Singular} helps to deduce $f = 0$ and the following
quadratic relations:
$$
\{ g, h, k \} \times \{ e, j \} = 0.
$$
If $(e,j) \neq (0,0)$, then $g = k = 0$ and so $e_{33}\not \in \Imm(R)$.
Hence, $e = j = 0$, and the following quadratic relations hold:
$$
\{ g, l, h+k \} \times (i-k+h) = 0, \quad
h(i+h) = gl.
$$
We have $g\neq0$, otherwise, the restrictions on $\Imm(R)$ do not work. Therefore, $k = i+h$ and $l = h(i+h)/g$.
Dividing by $g$, we get for $s = i/g$ and $t = h/g$
$$
R(e_{12}) = \begin{pmatrix}
s & 0 & 0 \\
0 & 0 & t(s+t) \\
0 & 0 & s+t
\end{pmatrix}, \quad
R(e_{13}) = \begin{pmatrix}
0 & 0 & 0 \\
0 & 0 & t \\
0 & 0 & 1
\end{pmatrix}, \quad s\neq0.
$$
Applying~$\psi$ with $\alpha = 1$, $\beta = \gamma = 0$, $\delta = 1/s$ and $\varepsilon = -t/s$, we get
\begin{equation}
R(e_{11}) = R(e_{22}) = R(e_{33}) = R(e_{23}) = 0, \quad
R(e_{12}) = e_{11} + e_{33}, \quad
R(e_{13}) = e_{33}.
\end{equation}

{\sc Subcase 2}: $h = k = l = 0$.
We compute that the following equations hold:
$$
\{ i,j \} \times \{ a,b,f,g \} = 0.
$$

{\sc Subcase 2.1}: $(i,j)\neq0$.
Then $a = b = f = g = 0$, and the four relations hold:
$$
\{ c, d, j \} \times e = 0, \quad
ic + jd = 0.
$$
We have $d\neq0$, otherwise $e_{33}\not\in\Imm(R)$.
Thus, $e = 0$ and $j = -ic/d$, i.\,e. after dividing by $d$ we get for $s = i/d$ and $t = c/d$
$$
R(e_{13}) = 0, \quad
R(e_{12}) = s(e_{11}-te_{13}),\quad
R(e_{23}) = e_{33} + te_{13}, \quad s\neq0.
$$
A conjugation with an automorphism~$\psi$ with nonzero coefficients $\alpha = \delta = 1/s$, $\gamma = -t$ gives us
\begin{equation}
R(e_{11}) = R(e_{22}) = R(e_{33}) = R(e_{13}) = 0, \quad
R(e_{12}) = e_{11}, \quad
R(e_{23}) = e_{33}.
\end{equation}

{\sc Subcase 2.2}: $i = j = 0$.
Then the quadratic relations include
$$
\{ a, e, f \} \times g = 0.
$$
If $g\neq0$, then $a = e = 0$ and $e_{11}\not\in\Imm(R)$.
Thus, $g = 0$, and the rest of the quadratic relations are
$$
\{ a, e, f \} \times c = 0, \quad
\{ b, e, a+f \} \times (f-a+d) = 0, \quad
f^2 + fd = be.
$$
Again, the conditions on~$\Imm(R)$ imply that $e\neq0$, thus, $c = 0$, $a = f + d$, and $b = f(f+d)/e$.
Dividing~$R$ by~$e$, we get for $s = f/e$ and $t = d/e$
$$
R(e_{12}) = 0, \quad
R(e_{13}) = e_{11} + se_{12}, \quad
R(e_{23}) = (s+t)(e_{11}+se_{12})+te_{33}, \quad t\neq0.
$$
Applying~$\psi$ with $\alpha = 1$, $\beta = s$, $\gamma = \varepsilon = 0$ and $\delta = t$, we get
\begin{equation}
R(e_{11}) = R(e_{22}) = R(e_{33}) = R(e_{12}) = 0, \quad
R(e_{13}) = e_{11}, \quad
R(e_{23}) = e_{11} + e_{33}.
\end{equation}

\section{Case $R(1) = e_{12}$}

We have $R(e_{12}) = 0$.
Computations with \texttt{Singular} give
\begin{gather*}
R(e_{11}) = \begin{pmatrix}
0 & 1-a-c & -b-d \\
0 & 0 & 0 \\
0 & 0 & 0
\end{pmatrix}, \quad
R(e_{22}) = \begin{pmatrix}
0 & a & b \\
0 & 0 & 0 \\
0 & 0 & 0
\end{pmatrix}, \quad
R(e_{33}) = \begin{pmatrix}
0 & c & d \\
0 & 0 & 0 \\
0 & 0 & 0
\end{pmatrix}, \\
R(e_{13}) = \begin{pmatrix}
0 & i-f & 0 \\
0 & 0 & 0 \\
0 & 0 & 0
\end{pmatrix}, \quad
R(e_{23}) = \begin{pmatrix}
f & g & h \\
0 & i & 0 \\
0 & 0 & j
\end{pmatrix}
\end{gather*}
joint with some quadratic relations, including
$$
\{f, j\} \times \{b, d, g\} = 0.
$$

{\sc Subcase 1}: $(f,j)\neq(0,0)$.
Then $b = d = g = 0$, and the rest of the quadratic relations are
$$
\{a,h\}\times (i-f) = 0, \quad
ij = 0, \quad
f(f-i-j) = 0,\quad
c(f-j) = 0.
$$

{\sc Subcase 1.1}: $j = 0$.
Then $f = i\neq 0$ and $c = 0$. It is an RB-operator
$R(e_{12}) = R(e_{13}) = R(e_{33}) = 0$, and
$$
R(e_{11}) = \begin{pmatrix}
0 & 1-a & 0 \\
0 & 0 & 0 \\
0 & 0 & 0
\end{pmatrix}, \quad
R(e_{22}) = \begin{pmatrix}
0 & a & 0 \\
0 & 0 & 0 \\
0 & 0 & 0
\end{pmatrix}, \quad
R(e_{23}) = \begin{pmatrix}
f & 0 & h \\
0 & f & 0 \\
0 & 0 & 0
\end{pmatrix}, \quad f\neq0.
$$
We apply conjugation with~$\psi\in\Aut(U_3(F))$
such that $\delta = 1$, $\alpha = \delta/f$, $\beta = \varepsilon = 0$ and $\gamma = h/f$.
If $a = 0$, we get
\begin{equation}
R(e_{12}) = R(e_{13}) = R(e_{22}) = R(e_{33}) = 0, \quad
R(e_{11}) = e_{12}, \quad
R(e_{23}) = e_{11}+e_{22}.
\end{equation}
If $a\neq0$, we take the same~$\psi$ with only one change: $\delta = a$ and get
\begin{equation}
R(e_{12}) {=} R(e_{13}) {=} R(e_{33}) {=} 0, \
R(e_{11}) {=} \lambda e_{12}, \
R(e_{22}) = e_{12}, \
R(e_{23}) = e_{11}+e_{22},\ \lambda\neq-1.
\end{equation}

{\sc Subcase 1.2}: $j \neq 0$.
Then $i = 0$, and the equations now have the form
$$
\{a,h,f-j\}\times f = 0, \quad
c(f-j) = 0.
$$

{\sc Subcase 1.2.1}: $f = 0$.
Then $c = 0$, and an RB-operator looks like $R(e_{12}) = R(e_{13}) = R(e_{33}) = 0$, and
$$
R(e_{11}) = \begin{pmatrix}
0 & 1-a & 0 \\
0 & 0 & 0 \\
0 & 0 & 0
\end{pmatrix}, \quad
R(e_{22}) = \begin{pmatrix}
0 & a & 0 \\
0 & 0 & 0 \\
0 & 0 & 0
\end{pmatrix}, \quad
R(e_{23}) = \begin{pmatrix}
0 & 0 & h \\
0 & 0 & 0 \\
0 & 0 & j
\end{pmatrix}, \quad j\neq0.
$$
Dividing by~$j$ and applying~$\psi$ with $\alpha = \delta$, $\beta = \varepsilon = 0$ and $\gamma = -h/j$, we get
\begin{equation*}
R(e_{12}) = R(e_{13}) = R(e_{33}) = 0, \
R(e_{11}) = x e_{12}, \
R(e_{22}) = y e_{12}, \
R(e_{23}) = e_{33}, \ y\neq-x.
\end{equation*}
If $y = 0$, we restrict the last $\psi$ with the condition $\delta = x$. Thus, we have
\begin{equation}
R(e_{12}) = R(e_{13}) = R(e_{22}) = R(e_{33}) = 0, \
R(e_{11}) = e_{12}, \quad
R(e_{23}) = e_{33}.
\end{equation}
Otherwise, we come to
\begin{equation}
R(e_{12}) = R(e_{13}) = R(e_{33}) = 0, \
R(e_{11}) = x e_{12}, \
R(e_{22}) = e_{12}, \
R(e_{23}) = e_{33}, \ x\neq-1.
\end{equation}

{\sc Subcase 1.2.2}: $f \neq 0$.
Then $f = j$, $a = h = 0$, and $R$ after dividing by~$f$ has the form that $R(e_{12}) = R(e_{22}) = 0$, and
$$
R(e_{11}) = (1-c)/f\, e_{12}, \quad
R(e_{33}) = (c/f)e_{12},\quad
R(e_{13}) = -e_{12},\quad
R(e_{23}) = e_{11}+e_{33}.
$$
Conjugate with~$\psi$, where
$\alpha = \delta = 1/f$, $\beta = \varepsilon = 0$, $\gamma = (1-c)/f$, and get
\begin{equation}
R(e_{12}) = R(e_{11}) = R(e_{22}) = 0 , \quad
R(e_{33}) = e_{12},\
R(e_{13}) = -e_{12},\
R(e_{23}) = e_{11}+e_{33}.
\end{equation}

{\sc Subcase 2}: $f = j = 0$.
Computations with \texttt{Singular} give
that $R(e_{12}) = 0$, and
\begin{gather*}
R(e_{11}) = \begin{pmatrix}
0 & 1-a-c & -b-d \\
0 & 0 & 0 \\
0 & 0 & 0
\end{pmatrix}, \quad
R(e_{22}) = \begin{pmatrix}
0 & a & b \\
0 & 0 & 0 \\
0 & 0 & 0
\end{pmatrix}, \quad
R(e_{33}) = \begin{pmatrix}
0 & c & d \\
0 & 0 & 0 \\
0 & 0 & 0
\end{pmatrix}, \\
R(e_{13}) = \begin{pmatrix}
0 & i & 0 \\
0 & 0 & 0 \\
0 & 0 & 0
\end{pmatrix}, \quad
R(e_{23}) = \begin{pmatrix}
0 & g & h \\
0 & i & 0 \\
0 & 0 & 0
\end{pmatrix}
\end{gather*}
joint with the following quadratic relations
$$
\{ a,b,d,g,h  \}\times i = 0.
$$

{\sc Subcase 2.1}: $i\neq 0$.
Then $a = b = d = g = h = 0$, and after deleting on $i$, we get that $R$ satisfies
$$
R(e_{12}) {=} R(e_{22}) = 0, \
R(e_{11}) = x e_{12}, \
R(e_{33}) = y e_{12}, \
R(e_{13}) = e_{12}, \
R(e_{23}) = e_{22}, \ x+y\neq0.
$$
Applying the automorphism $\psi$ with the parameters~$\gamma = (y-x)/2$, $\alpha = \delta = (x+y)/2$, $\beta = \varepsilon = 0$, we get the final form of~$R$:
\begin{equation}
R(e_{12}) = R(e_{22}) = 0, \quad
R(e_{11}) = e_{12}, \
R(e_{33}) = e_{12}, \
R(e_{13}) = e_{12}, \
R(e_{23}) = e_{22}.
\end{equation}

{\sc Subcase 2.2}: $i = 0$.
Then $R$ has the form
$$
R(e_{12}) = R(e_{13}) = 0, \quad
R(e_{11}), R(e_{22}), R(e_{33}), R(e_{23})\in L(e_{12},e_{13}).
$$

\section{Case $R(1) = e_{13}$}

Here $R(e_{13}) = 0$.
Computations with \texttt{Singular} imply
\begin{gather*}
R(e_{11}) = \begin{pmatrix}
-a-f & -b-g & 1-c-h \\
0 & -d-i & -e-j \\
0 & 0 & -a-f
\end{pmatrix}, \
R(e_{22}) = \begin{pmatrix}
a & b & c \\
0 & d & e \\
0 & 0 & a
\end{pmatrix}, \
R(e_{33}) = \begin{pmatrix}
f & g & h \\
0 & i & j \\
0 & 0 & f
\end{pmatrix}\!, \\
R(e_{12}) = \begin{pmatrix}
p & a-d+f-i & r \\
0 & s & t \\
0 & 0 & p
\end{pmatrix}, \quad
R(e_{23}) = \begin{pmatrix}
k & l & m \\
0 & n & i-f \\
0 & 0 & k
\end{pmatrix}
\end{gather*}
joint with some quadratic relations, including
$$
\{ f, i, g, k, l, n \} \times \{ a+f, d+i, e+j, p, r, s, t \} = 0.
$$

{\sc Subcase 1}: $(f,i,g,k,l,n)\neq(0,0,0,0,0,0)$.
Then $p = r = s = t = 0$, $a=-f$, $d=-i$, and $e = -j$.
Thus, $R(e_{12}) = R(e_{13}) = 0$,
\begin{gather*}
R(e_{11}) = \begin{pmatrix}
0 & -b-g & 1-c-h \\
0 & 0 & 0 \\
0 & 0 & 0
\end{pmatrix}, \quad
R(e_{22}) = \begin{pmatrix}
-f & b & c \\
0 & -i & -j \\
0 & 0 & -f
\end{pmatrix}, \\
R(e_{33}) = \begin{pmatrix}
f & g & h \\
0 & i & j \\
0 & 0 & f
\end{pmatrix}, \quad
R(e_{23}) = \begin{pmatrix}
k & l & m \\
0 & n & i-f \\
0 & 0 & k
\end{pmatrix},
\end{gather*}
and quadratic relations include
$$
\{ f, i, j, k, n \} \times \{ b+g, m-g\} = 0.
$$

{\sc Subcase 1.1}: $(f,i,j,k,n)\neq(0,0,0,0,0)$.
Hence, $m = -b = g$, and the relations include
$$
\{ f, k \} \times \{ i, n \} = 0.
$$

{\sc Subcase 1.1.1}: $(f,k)\neq(0,0)$.
Then $i = n = 0$, and the quadratic equations are
$$
kj + f^2 = 0, \quad
lj + fg = 0, \quad
gj - cf = 0, \quad
kg - lf = 0, \quad
ck + fg = 0.
$$
If $k = 0$, then $f = 0$, a contradiction.
So, $k\neq0$, and we express
$$
j = -f^2/k, \quad
g = lf/k, \quad
c = -f^2l/k^2,
$$
coming to the following form of~$R$ after dividing by~$k$:
($s = f/k$, $t = l/k$, \linebreak
$q = g/k$, $u = h/k$, $z = 1/k +s^2t -u$)
\begin{gather*}
R(e_{12}) = R(e_{13}) = 0, \quad
R(e_{11}) = \begin{pmatrix}
0 & 0 & z \\
0 & 0 & 0 \\
0 & 0 & 0
\end{pmatrix}, \quad
R(e_{22}) = \begin{pmatrix}
-s & -st & -s^2t \\
0 & 0 & s^2 \\
0 & 0 & -s
\end{pmatrix}, \\
R(e_{33}) = \begin{pmatrix}
s & st & u \\
0 & 0 & -s^2 \\
0 & 0 & s
\end{pmatrix}, \quad
R(e_{23}) = \begin{pmatrix}
1 & t & st \\
0 & 0 & -s \\
0 & 0 & 1
\end{pmatrix}.
\end{gather*}

Applying~$\psi$ with $\alpha = \delta = 1$,
$\beta = t$, $\gamma = t$, $\varepsilon = s$, we get for $z+y\neq0$
\begin{gather*}
R(e_{12}) = R(e_{13}) = R(e_{22}) = 0, \quad
R(e_{11}) = ze_{13}, \quad
R(e_{33}) = ye_{13}, \quad
R(e_{23}) = e_{11}+e_{33}.
\end{gather*}
Depending on whether $z$ or $y$ is nonzero, we get the two variants for $y\neq-1$:
\begin{gather}
R(e_{12}) = R(e_{13}) = R(e_{22}) = 0, \
R(e_{11}) = e_{13}, \
R(e_{33}) = ye_{13}, \
R(e_{23}) = e_{11}+e_{33}; \\
R(e_{12}) = R(e_{13}) = R(e_{22}) = 0, \
R(e_{11}) = ye_{13}, \
R(e_{33}) = e_{13}, \
R(e_{23}) = e_{11}+e_{33}.
\end{gather}

{\sc Subcase 1.1.2}: $f = k = 0$.
Then the quadratic relations are the following:
$$
i^2 = jn, \quad
gi = lj, \quad
ci = -gj, \quad
ng = il, \quad
cn = -jl.
$$

If $j = 0$, then $i = 0$.
The conditions of subcase~1.1 imply $n\neq0$, hence,
$c = g = 0$. We get an RB-operator
$$
R(e_{12}) = R(e_{13}) = R(e_{22}) = 0, \
R(e_{11}) = (1-h)e_{13}, \
R(e_{33}) = he_{13}, \
R(e_{23}) = le_{12} + ne_{22}.
$$
Let us split this variant into two.
If $h = 0$, then applying~$\psi$ with $\alpha = 1$, $\beta = -l/n$, $\gamma = \varepsilon = 0$ and $\delta = n$, we get
\begin{equation}\label{R(1)=e13,Subcase112}
R(e_{12}) = R(e_{13}) = R(e_{22}) = R(e_{33}) = 0, \quad
R(e_{11}) = e_{13}, \quad
R(e_{23}) = e_{22}.
\end{equation}
If $h \neq 0$, we apply~$\psi$ with $\alpha = h$, $\beta = -l/n$, $\gamma = \varepsilon = 0$, $\delta = nh$ and write down~$R$ 
depending on a parameter~$\varkappa\neq-1$:
\begin{equation}
R(e_{12}) = R(e_{13}) = R(e_{22}) = 0, \quad
R(e_{11}) = \varkappa e_{13}, \quad
R(e_{23}) = e_{22}, \quad
R(e_{33}) = e_{13}.
\end{equation}

If $j\neq0$, we express
$$
n = i^2/j, \quad
g = -ci/j, \quad
l = gi/j = -ci^2/j^2.
$$
Dividing by~$j$, we have an RB-operator~$R$ such that
$R(e_{12}) = R(e_{13}) = 0$ and
($s = i/j$, $t = c/j$, $q = h/j$)
\begin{gather*}
R(e_{11}) = \begin{pmatrix}
0 & 0 & 1/j-t-q \\
0 & 0 & 0 \\
0 & 0 & 0
\end{pmatrix}, \quad
R(e_{22}) = \begin{pmatrix}
0 & st & t \\
0 & -s & -1 \\
0 & 0 & 0
\end{pmatrix}, \\
R(e_{33}) = \begin{pmatrix}
0 & -st & q \\
0 & s & 1 \\
0 & 0 & 0
\end{pmatrix}, \quad
R(e_{23}) = \begin{pmatrix}
0 & -s^2t & -st \\
0 & s^2 & s \\
0 & 0 & 0
\end{pmatrix}.
\end{gather*}

When $s = 0$, we get
$$
R(e_{12}) = R(e_{13}) = R(e_{23}) = 0. \quad
R(e_{11}), R(e_{22}), R(e_{33})\in L(e_{13},e_{23}).
$$

Suppose that~$s\neq0$.
Applying~$\psi\in\Aut(U_3(F))$ with
$\alpha = 1/j$, $\beta = t$, $\gamma = 0$,
$\delta = s^2/j$,
$\varepsilon = s/j$,
we get
$$
R(e_{12}) {=} R(e_{13}) {=} R(e_{22}) = 0, \
R(e_{11}) = (1{-}jt{-}jq)e_{13}, \
R(e_{33}) = j(t{+}q)e_{13},\
R(e_{23}) = e_{22}.
$$
If $t+q = 0$, then we get again the RB-operator~\eqref{R(1)=e13,Subcase112}.
Otherwise, we apply~$\psi$ with $\alpha = j(t+q)$ and get~$R$ depending on a parameter~$\varkappa\neq-1$:
\begin{equation}
R(e_{12}) = R(e_{13}) = R(e_{22}) = 0, \quad
R(e_{11}) = \varkappa e_{13}, \quad
R(e_{33}) = e_{13},\quad
R(e_{23}) = e_{22}.
\end{equation}

{\sc Subcase 1.2}: $f = i = j = k = n = 0$.
It is an RB-operator~$R$ such that
$$
R(e_{12}) = R(e_{13}) = 0, \quad
R(e_{11}),
R(e_{22}),
R(e_{33}),
R(e_{23})\in L(e_{12},e_{13}).
$$

{\sc Subcase 2}: $f = i = g = k = l = n = 0$.
Hence, $R(e_{13}) = 0$, and
\begin{gather*}
R(e_{11}) = \begin{pmatrix}
-a & -b & 1-c-h \\
0 & -d & -e-j \\
0 & 0 & -a
\end{pmatrix}, \quad
R(e_{22}) = \begin{pmatrix}
a & b & c \\
0 & d & e \\
0 & 0 & a
\end{pmatrix}, \quad
R(e_{33}) = \begin{pmatrix}
0 & 0 & h \\
0 & 0 & j \\
0 & 0 & 0
\end{pmatrix}, \\
R(e_{12}) = \begin{pmatrix}
p & a-d & r \\
0 & s & t \\
0 & 0 & p
\end{pmatrix}, \quad
R(e_{23}) = \begin{pmatrix}
0 & 0 & m \\
0 & 0 & 0 \\
0 & 0 & 0
\end{pmatrix},
\end{gather*}
the quadratic relations include
$$
\{ j, m \} \times \{ a, d, p, s \} = 0.
$$

{\sc Subcase 2.1}: $(j,m) \neq (0,0)$.
Then $a = d = p = s = 0$, and the following equations hold:
$$
\{ e, j, t \} \times \{ b, m \} = 0,
\quad br = 0.
$$
If $(e,j,t)\neq(0,0,0)$, then $b = m = 0$, it is an RB-operator
\begin{gather*}
R(e_{13}) = R(e_{23}) = 0, \quad
R(e_{11}), R(e_{22}), R(e_{33}), R(e_{12}) \in L(e_{13},e_{23}).
\end{gather*}

Let $e = j = t = 0$.
If $b = 0$, then
\begin{gather*}
R(e_{13}) = 0, \quad
R(e_{11}),
R(e_{22}),
R(e_{33}),
R(e_{12}),
R(e_{23})\in L(e_{13}).
\end{gather*}

If $r = 0$, then
$$
R(e_{12}) = R(e_{13}) = 0, \quad
R(e_{11}),
R(e_{22}),
R(e_{33}),
R(e_{23}) \in L(e_{12},e_{13}).
$$

{\sc Subcase 2.2}: $j = m = 0$.
The following equations hold:
$$
\{ d,s \} \times \{ a,p,r+e\} = 0.
$$

{\sc Subcase 2.2.1}: $(d,s) \neq (0,0)$.
Then $a = p = 0$ and $r = -e$.
Additionally, we have the quadratic equations
$$
cd = be, \quad
dt = se, \quad
sc = -de, \quad
tb = -de, \quad
d^2 = -sb.
$$
Since $(d,s)\neq(0,0)$, we conclude that $s\neq0$.
Therefore,
$$
b = -d^2/s, \quad
c = -d^2t/s^2, \quad
e = dt/s.
$$
After dividing by~$s$, we get
$R(e_{13}) = R(e_{23}) = 0$, $R(e_{33}) = ze_{13}$,
\begin{gather*}
R(e_{11}) = \begin{pmatrix}
0 & -x^2 & q \\
0 & -x & -xy \\
0 & 0 & 0
\end{pmatrix}\!, \quad
R(e_{22}) = \begin{pmatrix}
0 & -x^2 & -x^2y \\
0 & x & xy \\
0 & 0 & 0
\end{pmatrix}\!, \quad
R(e_{12}) = \begin{pmatrix}
0 & -x & -xy \\
0 & 1 & y \\
0 & 0 & 0
\end{pmatrix}\!.
\end{gather*}
Conjugation with~$\psi$ with
$\alpha = \delta = 1$, $\beta = x$, $\gamma = 0$, $\varepsilon = y$ gives
$$
R(e_{13}) = R(e_{22}) = R(e_{23}) = 0,\quad R(e_{11}),R(e_{33})\in L(e_{13}), \quad
R(e_{12}) = e_{22}.
$$
However, $e_{11}\not\in\ker(P)$ or $e_{33}\not\in\ker(P)$.
Depending on the variant and involving the parameter~$\alpha$, we have one of the two possibilities with $x\neq-1$:
\begin{gather}
R(e_{13}) = R(e_{22}) = R(e_{23}) = 0,\quad
R(e_{11}) = e_{13}, \quad
R(e_{33}) = x e_{13}, \quad
R(e_{12}) = e_{22}; \\
R(e_{13}) = R(e_{22}) = R(e_{23}) = 0,\quad
R(e_{11}) = x e_{13}, \quad
R(e_{33}) = e_{13}, \quad
R(e_{12}) = e_{22}.
\end{gather}

{\sc Subcase 2.2.2}: $d = s = 0$.
Among the relations, we have
$$
\{ a,b,p \} \times (r+e) = 0.
$$

If $a = b = p = 0$, then it is an RB-operator
\begin{gather*}
R(e_{13}) = R(e_{23}) = 0, \quad
R(e_{11}), R(e_{22}), R(e_{33}), R(e_{12}) \in L(e_{13},e_{23}).
\end{gather*}

Suppose that $(a,b,p)\neq(0,0,0)$.
Hence, $r = -e$, and the list of quadratic relations consists of
$$
be+ca = 0, \quad
pe = ta, \quad
pc+ea = 0, \quad
a^2 = pb, \quad
tb = ea.
$$

{\sc Subcase 2.2.2.1}: $p = 0$.
Then $a = 0$. Since $(a,b,p)\neq(0,0,0)$, we conclude that $b = 0$ and so, $e = t = 0$.
We get an RB-operator~$R$ such that
$$
R(e_{12}) = R(e_{13}) = 0, \quad
R(e_{11}),
R(e_{22}),
R(e_{33}),
R(e_{23}) \in L(e_{12},e_{13}).
$$

{\sc Subcase 2.2.2.2}: $p \neq 0$.
Then we express
$$
b = a^2/p, \quad
c = -a^2t/p^2, \quad
e = at/p,
$$
it is~$R$ satisfying the conditions (after dividing by~$p$)
$R(e_{13}) = R(e_{23}) = 0$, $R(e_{33}) = ze_{13}$,
\begin{gather*}
R(e_{11}) = \begin{pmatrix}
-x & -x^2 & q \\
0 & 0 & -xy \\
0 & 0 & -x
\end{pmatrix}, \
R(e_{22}) = \begin{pmatrix}
x & x^2 & -x^2y \\
0 & 0 & xy \\
0 & 0 & x
\end{pmatrix}, \
R(e_{12}) = \begin{pmatrix}
1 & x & -xy \\
0 & 0 & y \\
0 & 0 & 1
\end{pmatrix}.
\end{gather*}
Conjugation with $\psi$ defined by~$\delta = 1$, $\beta = x$, $\gamma = 0$, $\varepsilon = -y$,
gives~$P$ such that
$e_{13},e_{22},e_{23}\in\ker(P)$,
$P(e_{11}),P(e_{33})\in L(e_{13})$,
$P(e_{12}) = e_{11}+e_{33}$.
However, $e_{11}\not\in\ker(P)$ or $e_{33}\not\in\ker(P)$.
Depending on the variant and involving the parameter~$\alpha$, we have one of the two possibilities with $x\neq-1$:
\begin{gather}
R(e_{13}) = R(e_{22}) = R(e_{23}) = 0, \
R(e_{11}) = e_{13}, \
R(e_{33}) = x e_{13}, \
R(e_{12}) = e_{11}+e_{33}; \\
R(e_{13}) = R(e_{22}) = R(e_{23}) = 0, \
R(e_{11}) = x e_{13}, \
R(e_{33}) = e_{13}, \
R(e_{12}) = e_{11}+e_{33}.
\end{gather}

\section{Case $R(1) = e_{12}+e_{23}$}

Fulfilling computations with \texttt{Singular},
we get among the linear relations the following one~$r_{11,12}=r_{33,23}-1$ and the unique quadratic equation
$r_{33,23}(r_{33,23}-1) = 0$.
Thus, up to the conjugation with~$\Theta_{1,3}$, we may assume that $r_{33,23} = 0$.
Hence, we get the following data: 
$e_{13},e_{23}\in \ker(R)$, $R(e_{12})=(1/2)e_{13}$ and
$$
R(e_{11}) = \begin{pmatrix}
0 & 1 & b \\
0 & 0 & 1/2 \\
0 & 0 & 0
\end{pmatrix}, \quad
R(e_{22}) = \begin{pmatrix}
0 & 0 & f \\
0 & 0 & 1/2 \\
0 & 0 & 0
\end{pmatrix}, \quad
R(e_{33}) = \begin{pmatrix}
0 & 0 & -b-f \\
0 & 0 & 0 \\
0 & 0 & 0
\end{pmatrix}.
$$

Conjugation with an automorphism $\varphi$ of $M_3(F)$, where
\begin{gather*}
\varphi(e_{ii})=e_{ii},\ i=1,2,3,\quad
\varphi(e_{12})=e_{12},\quad
\varphi(e_{13})=e_{13}/2,\quad
\varphi(e_{23})=e_{23}/2,\quad
\end{gather*}
maps $R$ to the RB-operator $P$ on $U_3(F)$ defined as follows,
\begin{equation}
\begin{gathered}
P(e_{13}) = P(e_{23}) = 0, \quad P(e_{12}) = e_{13}, \\
P(e_{11}) = \begin{pmatrix}
0 & 1 & x \\
0 & 0 & 1 \\
0 & 0 & 0
\end{pmatrix}, \quad
P(e_{22}) = \begin{pmatrix}
0 & 0 & y \\
0 & 0 & 1 \\
0 & 0 & 0
\end{pmatrix}, \quad
P(e_{33}) = \begin{pmatrix}
0 & 0 & -x-y \\
0 & 0 & 0 \\
0 & 0 & 0
\end{pmatrix}.
\end{gathered}
\end{equation}

\section{Classification}

Now, we gather all obtained operators jointly.

{\bf Theorem 2}.
All nonzero RB-operators of weight zero on $U_3(F)$ over a field $F$ of characteristic~0
up to conjugation with (anti)automorphisms $U_3(F)$ and multiplication on a nonzero scalar are the following: (below $\varkappa\neq-1$)

\newpage

\begin{gather*}
(R1)\quad
R(e_{12}) = R(e_{13}) = 0, \quad
R(e_{11}),
R(e_{22}),
R(e_{33}),
R(e_{23})\in L(e_{12},e_{13}); \\
(R2)\quad
R(e_{13}) = 0, \quad
R(e_{12}), R(e_{11}), R(e_{22}), R(e_{23}), R(e_{33}) \in L(e_{13}); \\
(R3)\quad
R(e_{11}) = R(e_{22}) = R(e_{23}) = R(e_{33}) = 0, \quad
R(e_{12}) = - e_{11}, \quad
R(e_{13}) = e_{23}; \\
(R4)\quad
R(e_{11}) = R(e_{13}) = R(e_{23}) = 0, \quad
R(e_{12}) =  e_{11}, \quad
R(e_{22}) = -e_{23}, \quad
R(e_{33}) =  e_{23}; \\
(R5)\quad
R(e_{11}) = R(e_{22}) = R(e_{33}) = R(e_{13}) =  R(e_{23}) = 0, \quad
R(e_{12}) = e_{11}; \\
(R6)\quad
R(e_{11}) = R(e_{22}) = R(e_{33}) = R(e_{12}) = R(e_{23}) = 0, \quad
R(e_{13}) =  e_{11}; \\
(R7)\quad
R(e_{12}) = R(e_{11}) = R(e_{22}) = R(e_{33}) = 0, \quad
R(e_{13}) = e_{12}, \quad
R(e_{23}) = e_{22}; \\
(R8)\quad
R(e_{12}) = R(e_{13}) = R(e_{11}) = R(e_{22}) = R(e_{33}) = 0, \quad
R(e_{23}) = e_{22}; \\
(R9)\quad
R(e_{12}) = R(e_{13}) = R(e_{22}) = 0, \quad
R(e_{11}) = -e_{13}, \quad
R(e_{22}) = e_{13},\quad
R(e_{23}) = e_{22}; \\
(R10)\quad
R(e_{12}) = R(e_{23}) = R(e_{11}) = R(e_{22}) = R(e_{33}) = 0, \quad
R(e_{13}) = e_{11} + e_{22}; \\
(R11)\quad
R(e_{12}) = R(e_{13}) = R(e_{11}) = R(e_{22}) = R(e_{33}) = 0, \quad
R(e_{23}) = e_{11}+e_{22}; \\
(R12)\quad
R(e_{12}) = R(e_{13}) = R(e_{33}) = 0, \
R(e_{23}) = e_{11}+e_{22}, \
R(e_{11}) = e_{12}, \
R(e_{22}) = -e_{12}; \\
(R13)\quad
R(e_{12}) {=} R(e_{33}) {=} 0, \
R(e_{13}) = e_{12}, \
R(e_{23}) = e_{11} {+} e_{22}, \
R(e_{11}) = e_{13}, \
R(e_{22}) {=} {-}e_{13};\\
(R14)\quad
R(e_{12}) = R(e_{11}) = R(e_{22}) = R(e_{33}) = 0, \quad
R(e_{13}) = e_{12}, \quad
R(e_{23}) = e_{11} + e_{22}; \\
(R15)\quad
R(e_{23}) = R(e_{11}) = R(e_{22}) = R(e_{33}) = 0, \quad
R(e_{13}) = -e_{23}, \quad
R(e_{12}) = e_{11}+e_{33}; \\
(R16)\quad
R(e_{23}) {=} R(e_{22}) {=} 0, \
R(e_{13}) {=} -e_{23}, \
R(e_{12}) {=} e_{11}{+}e_{33}, \
R(e_{11}) {=} -e_{23}, \
R(e_{33}) {=} e_{23}; \\
(R17)\quad
R(e_{23}) = R(e_{13}) = R(e_{11}) = R(e_{22}) = R(e_{33}) = 0, \quad
R(e_{12}) = e_{11}+e_{33}; \\
(R18)\quad
R(e_{23}) = R(e_{13}) = R(e_{22}) = 0, \
R(e_{12}) = e_{11}+e_{33}, \
R(e_{11}) = -e_{13},\
R(e_{33}) = e_{13}; \\
(R19)\quad
R(e_{11}) = R(e_{22}) = R(e_{33}) = R(e_{23}) = 0, \quad
R(e_{12}) = e_{22}, \quad
R(e_{13}) = e_{11} + e_{22}; \\
(R20)\quad
R(e_{11}) = R(e_{22}) = R(e_{33}) = R(e_{13}) = 0, \quad
R(e_{12}) = e_{11}, \quad
R(e_{23}) = e_{11} + e_{22}; \\
(R21)\quad
R(e_{11}) = R(e_{22}) = R(e_{33}) = R(e_{12}) = 0, \quad
R(e_{13}) = e_{11}, \quad
R(e_{23}) = e_{22}; \\
(R22)\quad
R(e_{11}) = R(e_{22}) = R(e_{33}) = R(e_{23}) = 0, \quad
R(e_{12}) = e_{11} + e_{33}, \quad
R(e_{13}) = e_{33}; \\
(R23)\quad
R(e_{11}) = R(e_{22}) = R(e_{33}) = R(e_{13}) = 0, \quad
R(e_{12}) = e_{11}, \quad
R(e_{23}) = e_{33}; \\
(R24)\quad
R(e_{11}) = R(e_{22}) = R(e_{33}) = R(e_{12}) = 0, \quad
R(e_{13}) = e_{11}, \quad
R(e_{23}) = e_{11} + e_{33}; \\
(R25)\quad
R(e_{12}) = R(e_{13}) = R(e_{22}) = R(e_{33}) = 0, \quad
R(e_{11}) = e_{12}, \quad
R(e_{23}) = e_{11}+e_{22}; \\
(R26)\quad
R(e_{12}) {=} R(e_{13}) {=} R(e_{33}) {=} 0, \
R(e_{11}) {=} \varkappa e_{12}, \
R(e_{22}) = e_{12}, \
R(e_{23}) = e_{11}+e_{22}; \\
(R27)\quad
R(e_{12}) = R(e_{13}) = R(e_{22}) = R(e_{33}) = 0, \
R(e_{11}) = e_{12}, \quad
R(e_{23}) = e_{33};\\
(R28)\quad R(e_{12}) = R(e_{13}) = R(e_{33}) = 0, \
R(e_{11}) = \varkappa e_{12}, \
R(e_{22}) = e_{12}, \
R(e_{23}) = e_{33};\\
(R29)\quad
R(e_{12}) {=} R(e_{11}) {=} R(e_{22}) = 0 , \quad
R(e_{33}) = e_{12},\
R(e_{13}) = {-}e_{12},\
R(e_{23}) = e_{11}{+}e_{33}; \\
(R30)\quad
R(e_{12}) = R(e_{22}) = 0, \
R(e_{11}) = e_{12}, \
R(e_{33}) = e_{12}, \
R(e_{13}) = e_{12}, \
R(e_{23}) = e_{22}; \\
(R31)\quad
R(e_{12}) = R(e_{13}) = R(e_{22}) = 0, \
R(e_{11}) = e_{13}, \
R(e_{33}) = \varkappa e_{13}, \
R(e_{23}) = e_{11}+e_{33}; \\
(R32)\quad
R(e_{12}) = R(e_{13}) = R(e_{22}) = 0, \
R(e_{11}) = \varkappa e_{13}, \
R(e_{33}) = e_{13}, \
R(e_{23}) = e_{11}+e_{33};\\
(R33)\quad
R(e_{12}) = R(e_{13}) = R(e_{22}) = R(e_{33}) = 0, \quad
R(e_{11}) = e_{13}, \quad
R(e_{23}) = e_{22};\\
(R34)\quad
R(e_{12}) = R(e_{13}) = R(e_{22}) = 0, \quad
R(e_{11}) = \varkappa e_{13}, \quad
R(e_{23}) = e_{22}, \quad
R(e_{33}) = e_{13};\\
(R35)\quad
R(e_{12}) = R(e_{13}) = R(e_{22}) = 0, \quad
R(e_{11}) = \varkappa e_{13}, \quad
R(e_{33}) = e_{13},\quad
R(e_{23}) = e_{22};\\
(R36)\quad
R(e_{13}) = R(e_{22}) = R(e_{23}) = 0,\quad
R(e_{11}) = e_{13}, \quad
R(e_{33}) = \varkappa e_{13}, \quad
R(e_{12}) = e_{22};  \\
(R37)\quad
R(e_{13}) = R(e_{22}) = R(e_{23}) = 0,\quad
R(e_{11}) = \varkappa e_{13}, \quad
R(e_{33}) = e_{13}, \quad
R(e_{12}) = e_{22}; \allowdisplaybreaks \\
(R38)\quad
R(e_{13}) = R(e_{22}) = R(e_{23}) = 0, \
R(e_{11}) = e_{13}, \
R(e_{33}) = \varkappa e_{13}, \
R(e_{12}) = e_{11}+e_{33}; \\
(R39)\quad
R(e_{13}) = R(e_{22}) = R(e_{23}) = 0, \
R(e_{11}) = \varkappa e_{13}, \
R(e_{33}) = e_{13}, \
R(e_{12}) = e_{11}+e_{33};\\
(R40)\quad
R(e_{13}) = R(e_{23}) = 0, \quad R(e_{12}) = e_{13}, \\
R(e_{11}) = \begin{pmatrix}
0 & 1 & b \\
0 & 0 & 1 \\
0 & 0 & 0
\end{pmatrix}, \quad
R(e_{22}) = \begin{pmatrix}
0 & 0 & f \\
0 & 0 & 1 \\
0 & 0 & 0
\end{pmatrix}, \quad
R(e_{33}) = \begin{pmatrix}
0 & 0 & -b-f \\
0 & 0 & 0 \\
0 & 0 & 0
\end{pmatrix}.
\end{gather*}

{\bf Corollary 1}.
All RB-operators $R$ of weight zero on $U_3(F)$, except (R40) (i.\,e. such that $R(1)^2\neq0$),
satisfy the property $\dim(\Imm(R))\leq 2$.

Let $A$ be an algebra.
Due to~\cite{Unital}, a Rota---Baxter index $\rb(A)$ of $A$ is defined as follows,
$$
\rb(A) = \min\{n\in\mathbb{N}\mid R^n = 0\mbox{ for an RB-operator }R \mbox{ of weight zero on }A\}.
$$
If such a number is undefined, put $\rb(A) = \infty$.

By Theorem~1, $\rb(U_3(F))<\infty$.
After completing the classification, we may compute this invariant directly.

{\bf Corollary 2}.
We have $\rb(U_3(F)) = 3$.

{\sc Proof}.
Note that $R^3 = 0$ holds for all operators from Theorem~3.
However, the RB-operators (R13), (R29), (R31), (R32), (R38), (R39), and (R40) satisfy $R^2\neq0$.
\hfill $\square$

\section*{Acknowledgements}

The study was supported by a grant from the Russian Science Foundation № 23-71-10005, https://rscf.ru/project/23-71-10005/

\noindent Vsevolod Gubarev \\
Novosibirsk State University \\
Pirogova str. 1, 630090 Novosibirsk, Russia \\
Sobolev Institute of Mathematics \\
Acad. Koptyug ave. 4, 630090 Novosibirsk, Russia \\
e-mail: wsewolod89@gmail.com

\end{document}